\documentclass{article}

\usepackage[OT1]{fontenc}
\usepackage{amsthm,amsmath}
\usepackage[numbers]{natbib}
\usepackage[colorlinks,citecolor=blue,urlcolor=blue]{hyperref}
\usepackage{amsfonts}
\usepackage[latin1]{inputenc}
\usepackage{amsbsy}
\usepackage{amssymb}
\usepackage{enumerate}
\usepackage{amsxtra}
\usepackage{amstext}
\usepackage{graphicx,graphics}
\usepackage{color}
\usepackage{epsfig}

\numberwithin{equation}{section}
\theoremstyle{plain}

\newtheorem{th-a}{Theorem}[subsection]
\newtheorem{lemma}{Lemma}
\newtheorem{lemma-a}{Lemma}[subsection]
\newtheorem{proposition}{Proposition}

\newcommand{\T}[1]{{\mathbb{T}}}

\def\Ccl{\mathcal{C}_{c\ell}}
\def\Ccld{\mathcal{C}_{c\ell}(\d)}
\def\leq{\leqslant}

\def\ni{\noindent}
\def\X{{\mathcal{X}}}

\def\a{\alpha}

\def\d{\delta}
\def\k{\kappa}

\def\1{1\!{\rm l}}

\def\ef{\tilde{f}}

\def\T{{\mathcal{T}}}

\def\argmin{{\mathop{{\mathrm{argmin}}}}}

\def\P{{{\mathrm{P}}}}
\def\E{{{\mathbb{E}}}}

\def\X{{\mathcal{X}}}

\begin{document}

\title{Risk Bounds for Embedded Variable Selection in Classification Trees}
\author{Servane Gey$^1$ \& Tristan Mary-Huard$^2$\\
Servane.Gey@parisdescartes.fr, maryhuar@agroparistech.fr\\
$^1$Laboratoire MAP5, UMR 8145, Universit\'{e} Paris Descartes, Paris, France\\
$^2$UMR AgroParisTech INRA MIA 518, Paris, France\\
}

\maketitle

\begin{abstract}
The problems of model and variable selections for classification trees are
jointly considered. A penalized criterion is proposed which explicitly
takes into account the number of variables, and a risk bound inequality is
provided for the tree classifier minimizing this criterion.
This penalized criterion is compared to the one used during the pruning step
of the CART algorithm. It is shown that the two criteria are
similar under some specific margin assumptions. In practice, the tuning
parameter of the CART penalty has to be calibrated by hold-out. Simulation
studies are performed which confirm that the hold-out procedure mimics the
form of the proposed penalized criterion.
\end{abstract}

\section{Introduction} \label{Sec : Intro}

Since the pioneering work of Breiman {\it et al.} \citep{Brei84}, classification trees have become a classical tool in machine learning. In particular, the Classification and Regression Tree (CART) algorithm is a
well-established algorithm to build and prune tree predictors. This algorithm
has been successfully applied in various fields, see for instance
\cite{Bel09,Chou89,DudFridSpee02,Wer98}.

\subsection{Building/selecting a tree}
The process of building (or choosing) a tree classifier from a training set can be summarized into an optimization problem, where the goal is to find the ``best'' tree classifier $\hat{f}$ satisfying
\begin{eqnarray}
\hat{f} = \arg\underset{f_T}{\min}\left(\P_nf_T + pen(n,T)\right) \ \ , \label{Equ:PenalizedCriterionMinimization}
\end{eqnarray}
where $n$ is the number of observations, $\P_n\hat{f}_T$ is the empirical risk of
tree classifier $\hat{f}_T$ based on tree $T$, and $pen(n,T)$ is a penalty function based on the size of the training set and on the characteristics of $T$.

Obtaining the best tree classifier $\hat{f}$ necessitates to solve a non-convex function over a large set of trees, something unfeasible in practice. As an alternative, a 2-step heuristic approach to solve this problem has been proposed in \citep{Brei84}, in the particular case where the penalized criterion is of the form
\begin{eqnarray}
pen(n,T) = \alpha_n\times |T| \ \ , \label{Eq : FormePenaliteCART}
\end{eqnarray}
where $\alpha_n$ is a tuning parameter that depends
on $n$, and $|T|$ is the size of the tree, i.e. the number of leaves (terminal nodes) of $T$. In the first step (called the \emph{growing step}) a large tree $T_{max}$ that achieves a perfect classification on the training set is built. Then, during the second step (called the \emph{pruning step}), the optimal subtree is obtained from the large tree, where the optimal subtree satisfies \begin{eqnarray*}
\hat{f}_{prun} = \underset{f_T, \ T\subseteq T_{max}}{\arg\min}\P_nf_T + \alpha_n\times |T| \ \ .
\end{eqnarray*}
While this heuristic approach is at the heart of the CART algorithm and is probably the most popular strategy to prune a tree, one should keep in mind that the actual goal is in fact to solve Problem \eqref{Equ:PenalizedCriterionMinimization}, and to obtain the properties of $\hat{f}$, whatever the (approximate) strategy that is applied to find it. \\

From a theoretical point of view, many works
have investigated the performance of the tree classifier resulting from the pruning step of CART rather than from the generic optimization problem. In the Gaussian or bounded regression context, penalty \eqref{Eq : FormePenaliteCART} was validated in
\cite{GeyNed05} using model selection framework. Another validation was obtained in the classification framework in \cite{Nob02}. More recently, a refined analysis of the pruning step was proposed in \cite{Gey12}, where margin adaptive risk bounds were obtained in the binary classification context. Importantly, these theoretical results are actually obtained conditionally to the construction of $T_{max}$. This means that only the performance of the pruning step is assessed, while the growing step is not taken into account.\\

\subsection{Classification trees and variable selection}
Because they are based on the 2-step heuristic of the CART algorithm, results obtained so far fail to take into account the complete process of obtaining a tree classifier. In particular, the embedded variable selection process that is inherent to tree classification algorithms has never been investigated. A variable selection
process is called embedded when it is included in the training step of the
classification algorithm. Therefore the learning and variable selection
processes cannot be separated. This embedded property is actually one of
the main arguments for the use of tree classifiers to deal with large dimension data (see
\cite{Brei01,DyeKahLeb07,GeyLeb08} for example). Note that in the CART algorithm, the inner selection process results
from the recursive growing strategy of the tree: at each node, the ``best''
variable is selected among all for splitting. As a result, in many cases the
maximal tree (and consequently all of its subtrees) only includes a small
subset of the $p$ initial variables. As a consequence, as long as tree classifiers are studied through the pruning step of the CART heuristic (hence conditionally to the growing step), it is impossible to investigate the complete variable selection process.\\

Although the embedded variable selection process is well-known (\cite{CzeWuWal08,GuyEli03,LalChaWesEli06}), it may appear at first glance that it is not correctly handled in the optimization program
\begin{eqnarray}
\hat{f} = \arg\underset{f_T}{\min}\left(\P_nf_T + \alpha_n\times |T| \right)\ \ , \label{Equ:CARTCriterionMinimization}
\end{eqnarray}
assuming the form of the penalty proposed in  \cite{Brei84} is correct. Indeed, this penalized criterion does not obviously depend on the total number of covariates $p$. This can be astonishing: in both the regression and
classification frameworks, theoretical studies have shown that in the variable
selection context, an extra term should be added to the penalty that is used
when only one model is considered per dimension (\cite{BirMas01,MarRobDau07})
to obtain oracle-type inequalities. Since the collection of possible trees increases with $p$, $p$ should play a crucial role in the regularization term.

Since parameter $p$ does not explicitly appear in criterion \eqref{Equ:CARTCriterionMinimization}, one can argue that $p$ is hidden in the constant term $\alpha_n$. This argument is verified from at least two penalties that can be exhibited from previous works:
\begin{itemize}
\item In \cite{Nob02} (equation 4), the penalty term has the form
\begin{eqnarray*}
pen(|T|,n) & = & C_1 \times \sqrt{\frac{|T|p\log n}{n}} \\
              & = & \hspace{-0.2cm}\sqrt{Cst\frac{p\log n}{n}} \times \sqrt{|T|}\\
              & = & \hspace{-0.2cm}\alpha(p,n)\sqrt{|T|}\ \ ,
\end{eqnarray*}
\item In \cite{Gey12} (Theorem 1), the penalty term is of order 
\begin{eqnarray*}
pen(|T|,n) & \approx & C_2 \times \frac{p\log{(p)}(1+\log(n/\log{(p)}))}{n} \times |T|\\
              & = & \hspace{-0.2cm}\alpha(p,n)|T|\ \ ,
\end{eqnarray*}
\end{itemize}
where $C_1$ and $C_2$ are known constants. While these two penalty functions depend on $p$, one can observe that their scaling order is much larger than the $\log(p)$ usually obtained in the variable selection context \cite{BirMas01,MarRobDau07}.\\

\subsection{Contribution}

The goal of the present paper is to investigate the classification performance of the tree classifier obtained by solving Problem \ref{Equ:PenalizedCriterionMinimization}, and to decipher the exact impact of variable selection on tree classifier selection. While this impact is theoretically studied through an ideal exhaustive selection procedure (unfeasible in practice), it sheds light on the heuristic procedures currently used in practice to mimic the ideal one (see Section \ref{Sec : simulations}). From a theoretical point of view, we consider the model selection problem where the goal is to select a candidate from \emph{all} possible tree classifiers. The strategy consists in choosing the candidate minimizing a penalized criterion that depends on parameters $p$ and $n$.  In this model selection context, we exhibit a penalization function where the variable selection process is explicitly taken into account, and provide performance guarantees for the candidate tree classifier through an upper bound of its risk. Then it is shown that the impact of variable selection, although investigated via the theoretical minimization problem \eqref{Equ:PenalizedCriterionMinimization}, can also be exhibited in practice for practical heuristic approaches. More precisely, a simulation study is performed which shows that the proposed theoretical penalization function is actually the one that is implicitly used in the pruning step of the CART algorithm.\\

The paper is organized as follows. Section \ref{Sec : Context} presents the framework of binary classification and describes tree classifiers. The main theoretical contribution and the simulation study are presented in Section \ref{Sec : RiskBounds}.
 Some discussion is developed in Section \ref{Sec : discussion},
and finally Section \ref{Sec : proofs} gives the proofs of the results
presented in Section \ref{Sec : RiskBounds}.

\section{Context} \label{Sec : Context}

\subsection{Classification framework}
The considered classification framework is the following. Suppose one observes a sample $\{(X_1,Y_1),\ldots,(X_n,Y_n)\}$ of $n$ independent copies of the random variable $(X,Y)$, where the
explanatory variable $X$ takes values in a measurable space $\X$ of dimension $p\geqslant 2$, and is
associated with a label $Y$ taking values in $\{0,1\}$. Suppose moreover that
each coordinate of $\X$ is ordered (i.e. $\X$ is a product of $p$ ordered subspaces). A  classifier is then any function $f$
mapping $\X$ into $\{0,1\}$. \\
The quality of a classifier is measured by its misclassification
rate
\begin{eqnarray} \label{Equ : misclassificationrate}
\P f := \P(f(X)\neq Y) \ \ ,
\end{eqnarray}
where $\P$ denotes the joint distribution of $(X,Y)$. If the joint
distribution of $(X,Y)$ were known, the problem
of finding an optimal classifier minimizing the misclassification rate would
be easily solved by considering the Bayes classifier $f^*$ defined for every
$x\in \X$ by
\begin{equation} \label{bayes}
f^*(x)=\1_{\eta(x)\geqslant 1/2} \ \ ,
\end{equation}
where $\eta(x)$ is the
conditional expectation of $Y$ given $X=x$, that is
\begin{equation} \label{eta}
\eta(x)=\P\left[Y=1 \ | \ X=x\right] \ \ .
\end{equation}
As $\P$ is unknown, the goal is to construct from sample
$\{(X_1,Y_1),\ldots,(X_n,Y_n)\}$ a classifier $\ef$ that is as close as
possible to $f^*$ in the following sense: since $f^*$ minimizes the
misclassification rate, $\ef$ will be chosen in such a way that its
misclassification rate is as close as possible to the misclassification rate
of $f^*$, i.e. in such a way that the loss
\begin{equation} \label{loss}
l(f^*,\ef)=\P(\ef(X)\neq Y)-\P(f^*(X)\neq Y)
\end{equation}
is as small as possible. \\
Many strategies or classification algorithms have been proposed to build $\ef$
(see  \cite{FrieHasTib01}, \cite{BouBouLug05} for an overview). The quality of a strategy is measured by its risk
\begin{eqnarray*}
\E[l(f^*,\ef)] \ \ ,
\end{eqnarray*}
where the expectation is taken with respect to the sample distribution. In the model selection framework, two strategies are usually considered:
\begin{itemize}
\item Empirical Risk Minimization: $\ef$ is chosen as the minimizer of
\begin{eqnarray} \label{Equ : empiricalrisk}
\P_n f := \frac{1}{n}\sum_{i=1}^n\1_{\{f(X_i)\neq Y_i\}} \ \ ,
\end{eqnarray}
over all classifiers $f$ belonging to a single class of classifiers,
\item Structural Risk Minimization: $\ef$ is chosen as the minimizer of the penalized empirical risk over a collection of classes.
\end{itemize}

\subsection{Margin assumptions}
It is now well known that without any assumption on the joint distribution $\P$, when considering a class of classifiers with finite Vapnik Chervonenkis (VC) dimension, the minimax convergence rate of the risk bound is of order $\mathcal{O}(1/\sqrt{n})$. It has also been shown that,
under the overoptimistic zero-error assumption (that is $Y=\eta(X)$ almost
surely, where $\eta$ is defined by \eqref{eta}), this minimax convergence rate
is at best of order $\mathcal{O}(1/n)$ (see \cite{Vap98, Lug02} for example).

These two extreme
cases can be modulated by so-called {\it margin assumptions} that make the link between the ``global'' pessimistic case (without any
assumption on $\P$) and the zero-error case (\cite{Kol06, Kol06rej, MamTsy99, MasNed06, Mas07, Tsy04, TsyGee05}).  \\

\ni In this paper, we consider the margin assumption proposed in \cite{MamTsy99}:
\begin{description}
\item {\bf MA(1)} There exist some constants $C_0>0$ and $\kappa>1$ such
  that, for all $t>0$,
\begin{eqnarray} \label{marginTsy}
\P\left(|2\eta(X)-1|\leqslant t\right)\leqslant C_0 \ t^{\frac{1}{\kappa-1}},
\end{eqnarray}
\end{description}
Note that by taking $t=h\in ]0,1[$ and the limit value
$\kappa=1$, we obtain the stronger assumption proposed in \cite{MasNed06} (see also the slightly weaker condition proposed in \cite{KohKrz07}):
\begin{description}
\item {\bf MA(2)} There exists $h\in ]0;1[$ such that
\begin{eqnarray} \label{marginMas}
\P\left(|2\eta(X)-1|\leq h\right)=0.
\end{eqnarray}
\end{description}
Assumption {\bf MA(2)} has an intuitive interpretation. It means that $(X, Y )$ is sufficiently well distributed to ensure that there
is no region in $X$ for which the toss-up strategy could be favored over others: $h$ can be
viewed as a measurement of the gap between labels 0 and 1 in the sense that, if $\eta(x)$ is
too close to 1/2, then choosing 0 or 1 will not make a real difference for that $x$. From a general point of view, the margin parameter quantifies the noise level of the classification problem, and may be understood as the equivalent of the variance parameter in the Gaussian model selection setting.

\subsection{Tree classifiers, classes of tree classifiers} \label{subsec : treeclass}
\begin{figure}[h]
\begin{center}
\includegraphics[scale = 0.25]{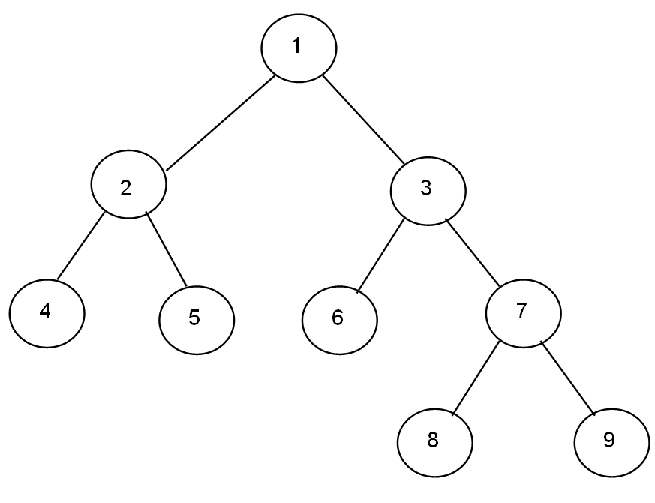}
\caption{Tree configuration example: for each node, the parent and child nodes are known. \label{ExempleCART} }
\end{center}
\end{figure}

\begin{figure}
\begin{center}
\includegraphics[scale=0.4]{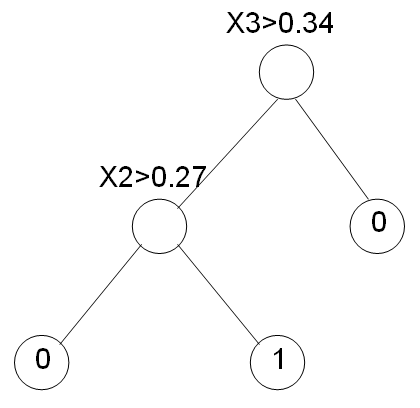}
\includegraphics[scale=0.4]{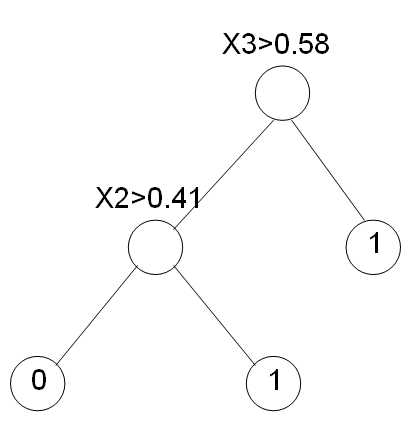}
\end{center}
\caption{Two tree classifiers that belong to the same class.}  \label{tree}
\end{figure}

A tree $T$ is a structure that can be represented as a hierarchy
whose elements are called nodes. For binary trees, each node has
either 0 or 2 children (called Left and Right). The initial node is
called the root of the tree and a node with no child is called a
leaf. The size of tree $T$ is defined as the number of its leaves
and noted $|T|$ in the following. In this paper, we define a tree $T_{c\ell}$ by two elements:
\begin{itemize}
\item its configuration $c$, i.e. the hierarchy between the nodes: for instance, in Figure \ref{ExempleCART}, we know that node 6 is the Left child node of node 3, and so on,
\item  the ordered list $\ell$ of variables that appear at each node, i.e. the $k^{th}$ variable in the list appears in node $k$.
\end{itemize}
A tree classifier $f$ based on tree $T_{c\ell}$ associates
\begin{itemize}
\item at each internal node a condition of the form $"X^{j_k}>s^k"$, where $j_{k}$ is the index of the variable associated with node $k$ and $s^k$ is a threshold,
\item  at each terminal node a label (here 0 or 1).
\end{itemize}
Therefore, an observation $x \in \mathcal{X}$ will be classified as follows: starting at the root, observation $x$
will move from a node of $f$ to another using the following rule: at
node $k$, if $"x^{j_k}>s^k"$ then $x$ moves to Right, otherwise it moves to
Left. At the end of the process, $x$ will be classified according to the label
of the leaf it reaches. \\
To summarize, a tree classifier associated with tree $T_{c\ell}$ splits $\mathcal{X}$ into $|T_{c\ell}|$ regions each associated with a label, and two classifiers associated with the same tree $T_{c\ell}$ differ in that the thresholds (for internal nodes) and labels (for leaves) are not the same. An example of two such tree classifiers is given in
Figure \ref{tree}.
\ni In the following, we will consider classes $\mathcal{C}_{c\ell} = \left\{f \  /  f \text{ based on } T_{c\ell}\right\}$ of classifiers based on a same tree $T_{c\ell}$.

\ni Finally, we define
\begin{eqnarray}
\overline{f}_{c\ell} \in \arg\underset{f\in\Ccl}{\min} \P f,
\end{eqnarray}
where $\P f$ is defined by \eqref{Equ : misclassificationrate}.

\section{Results} \label{Sec : RiskBounds}

\subsection{Risk bounds} \label{subsec : riskbound1class}

We first consider a single class $\mathcal{C}_{c\ell}$ of tree classifiers and
its associated empirical risk minimizer $$\widehat{f}_{c\ell} \in
\arg\underset{f\in\Ccl}{\min} \P_n f,$$ where $\P_n f$ is defined by \eqref{Equ : empiricalrisk}.

\begin{proposition}\label{Proposition : Control1class}
Assume that margin assumption {\bf MA(1)} is verified. For all $t_{c\ell}>0$
and $\alpha>0$, there exist positive constants  $K_1$, $K_2$, $K$ depending on $\alpha$, $C_0$ and $\kappa$ such that, with probability at least $1-e^{-t_{c\ell}}$,
\begin{eqnarray}
l(f^*,\widehat{f}_{c\ell}) \leq (1+\alpha)l(f^*,\overline{f}_{c\ell}) + K_1\left(\frac{|T_{c\ell}|\log(2n)}{n}\right)^\frac{\kappa}{2\kappa-1} + K_2\left(\frac{t_{c\ell}}{n}\right)^\frac{\kappa}{2\kappa-1} + K\frac{t_{c\ell}}{n} \ \ . \label{PropControl1Class : Traject}
\end{eqnarray}
Moreover, we obtain the following upper bound
\begin{eqnarray}
E\left[l(f^*,\widehat{f}_{c\ell})\right] \leq (1+\alpha)l(f^*,\overline{f}_{c\ell}) + K_1\left(\frac{|T_{c\ell}|\log(2n)}{n}\right)^\frac{\kappa}{2\kappa-1} + Cn^{-\frac{\kappa}{2\kappa-1}} \ \ . \label{PropControl1Class : Esperance}
\end{eqnarray}
\end{proposition}

The obtained bound is in keeping with classical results already given in
\cite{MamTsy99}. In particular, if the Bayes classifier belongs to class $\mathcal{C}_{c\ell}$, the rate of convergence for the risk associated with estimator $\widehat{f}_{c\ell}$ is of order $(\log(2n)/n)^\frac{\kappa}{2\kappa-1}$. \\

In practice, since no information is available about how to choose class
$\mathcal{C}_{c\ell}$, one needs to consider the collection $\mathcal{M}$ of
all possible configurations and variable lists. In each class
$\mathcal{C}_{c\ell}$, a candidate $\widehat{f}_{c\ell}$ is chosen by empirical risk minimization,
then the final classifier $\ef$ is selected among all class candidates by
minimization of a penalized criterion:
\begin{eqnarray*}
\widehat{c\ell} &=& \arg\underset{c,\ell}{\min}\left( \P_n\widehat{f}_{c\ell} +
pen(c,\ell)\right) \ \ ,\\
\ef &=&\widehat{f}_{\widehat{c\ell}} \ \ .
\end{eqnarray*}
The following result provides insight about how the penalty should be chosen to
ensure good performance for $\ef$.

\begin{proposition} \label{Proposition : ModelSelection}
Assume that margin assumption {\bf MA(1)} is
verified. If
\begin{eqnarray} \label{Eq : critere}
\ef=\underset{\{\widehat{f}_{c\ell} \ , \ (c,\ell)\in \mathcal{M}\}}{\argmin}\left(
  \P_n\widehat{f}_{c\ell} + pen(c,\ell)\right),
\end{eqnarray}
where
\begin{eqnarray} \label{Eq : penalite}
pen(c,\ell) &\hspace{-0.3cm}=& \hspace{-0.3cm} C'_\kappa\left(\frac{|T_{c\ell}|\log(2n)}{n}\right)^\frac{\kappa}{2\kappa-1} + C''_\kappa\left(\frac{|T_{c\ell}|\log(p)}{n}\right)^\frac{\kappa}{2\kappa-1}
\end{eqnarray}
with constants $C'_\kappa$ and $C''_\kappa$ depending on $C_0$ and $\k$
appearing in the margin condition, then there exist positive constants $C'_1, \ C'_2$ and $\Sigma$ such that with probability at least $1-3\Sigma e^{-x}$
\begin{eqnarray*}
l(f^*,\tilde{f}) \leq C_1'\underset{c,\ell}{\inf}\left\{\inf_{f\in\mathcal{C}_{c\ell}}l(f^*,f) +
pen(c,\ell) \right\}+ C'_2\left( \left(\frac{x}{n}\right)^\frac{\kappa}{2\kappa-1} +
\frac{x}{n}\right) \ \ .
\end{eqnarray*}
Moreover, we obtain the following upper bound:
\begin{eqnarray}
\E[l(f^*,\ef)]\leq  C'_1\inf_{(c,\ell)\in\mathcal{M}}\left\{\inf_{f\in\mathcal{C}_{c\ell}}l(f^*,f)+pen(c,\ell)\right\}+\frac{C''_2\times \Sigma}{n^{\frac{\k}{2\k-1}}}. \label{Eq : BorneRisqueProp2}
\end{eqnarray}
\end{proposition}
The proofs of Propositions \ref{Proposition : Control1class} and
\ref{Proposition : ModelSelection} are given in Section \ref{Sec : proofs}.\\

Several comments can be made about the result of Proposition \ref{Proposition : ModelSelection}:

\paragraph{Quality of the upper bound} Compared with previous results \cite{Nob02,Gey12}, the upper bound for the risk is improved in two different ways. First, since all possible binary trees are considered, in the present result the complete construction path of the tree classifier is taken into account: the infimum in equation \eqref{Eq : BorneRisqueProp2} is taken on all possible classes of tree classifiers. Conversely, in previous results only the performance of the pruning step was assessed, i.e. the corresponding infimum was restricted to the list of classes associated with subtrees of the maximal tree. Second, thanks to the margin hypothesis, the convergence rate of the upper bound is faster than $\mathcal{O}(1/\sqrt{n})$ as soon as $\kappa < +\infty$.


\paragraph{Margin parameter} The proposed penalty \eqref{Eq : penalite} depends on the margin parameter $\kappa$, that is usually unknown in practice. From a theoretical point of view, because this parameter quantifies the noise level of the classification problem, it necessarily appears in the ideal penalty function (as does the unknown variance in Gaussian model selection). From a practical point of view, it has to be estimated from the data. Obtaining this estimate in the general case is an open question.

\paragraph{Strong margin assumption} In the particular case of margin assumption
{\bf MA(2)} given by equation \eqref{marginMas}, penalty \eqref{Eq : penalite}
becomes (taking $\kappa=1$):
\begin{eqnarray*}
pen(c,\ell) & = & \frac{C^1_h\log(2n) + C^2_h\log(p)}{n}|T_{c\ell}|\\
            & = & \a_n|T_{c\ell}|.
\end{eqnarray*}
This corresponds exactly to the penalty proposed in \cite{Brei84} for the CART
algorithm (see equation \eqref{Eq : FormePenaliteCART}). This penalty function has already been validated for the pruning step of the CART algorithm, (see \cite{GeyNed05} for the regression framework  and \cite{Gey12} for the binary classification framework). A similar result is established by Proposition \ref{Proposition : ModelSelection} when considering the exact optimization problem \eqref{Equ:PenalizedCriterionMinimization}.  Also note that in this context the margin parameter only appears in constant $a_n$. Because this constant will be tuned accordingly to the data (using cross-validation for instance), the problem of estimating the margin parameter is discarded.


\paragraph{ Variable selection} In comparison with the upper bound  obtained in Proposition \ref{Proposition : Control1class}, one can observe in \eqref{Eq : BorneRisqueProp2} the impact of parameter $p$ that appears through the penalty. This quantity arises during the union bound step of the proof (see Section \ref{subsec : proofmodelselection}), where one has to count the number of classes sharing the same complexity. This conveys the fact that to build an optimal tree of size $k$, one has to choose $k$ variables among $p$ (with replacement). This is obviously a much easier task when $p=100$ than when $p=10,000$. This is where the variable selection task is taken into account. Moreover, the penalty term can be upper
bounded by
\begin{eqnarray*}
pen(c,\ell) \leq C'_\kappa\left(\frac{|T_{c\ell}|\log(2n)}{n}\right)^\frac{\kappa}{2\kappa-1} + \log(p)\left(C''_\kappa\left(\frac{|T_{c\ell}|}{n}\right)^\frac{\kappa}{2\kappa-1} +
C'''_\kappa\frac{|T_{c\ell}|}{n}\right),
\end{eqnarray*}
advocating for a penalty that should be linear with respect to $\log(p)$. This
linear relationship is investigated in Section \ref{Sec : simulations}.

\paragraph{ Oracle-type inequality} Vapnik-Chervonenkis bounds for binary classification
without any margin assumption give the following penalty form (see \cite{DevGyoLug96} for instance)
\begin{eqnarray*}
pen_{V}(c,\ell) = C^1_V\sqrt{\frac{|T_{c\ell}|\log(n)}{n}} + C^2_V\frac{|T_{c\ell}|}{n}.
\end{eqnarray*}
This implies that, for classes
associated with trees of large size, $pen(c,\ell)$ given in \eqref{Eq : penalite} becomes larger than
$pen_{V}(c,\ell)$. Therefore, to obtain an oracle-type inequality,
$pen(c,\ell)$ can be replaced by $\min\left\{pen_{V}(c,\ell),pen(c,\ell)\right\}$.

\subsection{Illustration on simulated data} \label{Sec : simulations}

\subsubsection{Practical determination of $\ef$}

The application of the strategy described in Proposition \ref{Proposition :
  ModelSelection} necessitates finding the empirical risk minimizer in each
class $\mathcal{C}_{c\ell}$, and then comparing all the candidates
$\widehat{f}_{c\ell}$ using the penalized criterion given by
\eqref{Eq : critere}. From a computational point of view, the exhaustive comparison
among all classes is an NP-hard problem. Therefore we need heuristic algorithms
to obtain a sequence of near-optimal penalized risk minimizers
$\left(\widehat{f}_k\right)_{k\geqslant 1}$ such that
\begin{eqnarray*}
\widehat{f}_{k} \approx \underset{\{ \widehat{f}_{c\ell}, \ |T_{c\ell}|=k\}}{\argmin} P_n \widehat{f}_{c\ell} \ \ .
\end{eqnarray*}
The CART algorithm, when applied with the empirical risk as an impurity
measure at each node (see \cite{FrieHasTib01}), may be understood as a forward
heuristic algorithm to build the sequence of optimal tree classifiers. In
particular, the subtree classifier $\widehat{f}_{k}$ of size $k$ extracted
from the maximal tree can be interpreted as the (approximate) optimizer of the
empirical risk over all the possible trees of size $k$. \\

This new understanding of the CART algorithm as a heuristic approach to obtain
the sequence of subtree minimizers is important, because it points out that
these subtree classifiers $\widehat{f}_{k}$ should be penalized as if the
exhaustive search were performed, {\it i.e.} using penalty given by \eqref{Eq
  : penalite}.\\

In most applications, when dealing with the construction of a tree classifier,
experimenters use criterion \eqref{Eq : FormePenaliteCART} in a growing-pruning strategy, and the unknown parameter $\alpha_n$ is chosen by hold-out or Q-fold cross-validation. This estimated value can be compared
with its theoretical counterpart given in \eqref{Eq : penalite}. To this end,
we perform a simulation study and compare the $\alpha_n$ obtained by
cross-validation to its theoretical form $$\frac{C^1_h\log(2n) +
  C^2_h\log(p)}{n}$$ obtained under the strong margin assumption {\bf MA(2)}.

\subsubsection{Simulations}

We consider four simulation designs:

\paragraph{Design 1} Variables $X^1,...,X^p$ are independently generated with
distribution $\mathcal{N}(0,1)$. The label is generated
as follows: If $X^1 > 0$ and $X^2>0$ then $Y=1$ with probability $q$, otherwise
$Y=1$ with probability $1-q$. Therefore only variables $X^1$ and $X^2$ are
informative. In this design, the Bayes classifier can be represented as a tree with 3 leaves, hence it belongs to the considered collection of classes. Moreover, variables are independent, and margin assumption {\bf MA(2)} is satisfied.

\paragraph{Design 2} First the labels are generated according to a
Bernoulli distribution with parameter $1/2$. Then variable $X^1$ is generated
such that $X^1|Y=0$  and  $X^1|Y=1$ are normally distributed with means $0$
and $1$, respectively, and variance $\sigma^2$. Variables $X^2,...,X^p$ are independent with
distribution $\mathcal{N}(0,1)$ and are non-informative. As for design 1,
the Bayes classifier can be represented as a tree and variables are independent, but it is
easy to show that margin assumption {\bf MA(2)} is not satisfied.

\paragraph{Design 3} Labels are simulated as in design 2. Then
variables $X^1$ and $X^2$ are generated such that, for $j=1,2$, $X^j|Y=0$  and  $X^j|Y=1$ are normally distributed with means $0$
and $1$, respectively, and variance $\sigma^2$. The last $p-2$ variables are
independent and non-informative. Here the Bayes classifier no longer belongs to
the collection of tree classes, and margin assumption {\bf MA(2)} is not satisfied.\\

\paragraph{Design 4} Three independent variables $X^1,X^2,X^3$ are generated
with distribution $\mathcal{N}(0,1)$. Each additional variable $X^j$ is then
simulated as a noisy copy of $(X^1+X^2+X^3)/\sqrt{3}$. The label is generated
as follows: If $(X^1)^2 + (X^2)^2 + (X^3)^2>2.5$ then $Y=1$, else $Y=0$. Here,
all the variables are correlated (with a strong correlation between the extra
variables), the Bayes classifier cannot be represented as a tree, and margin assumption {\bf
MA(2)} is not satisfied.\\

For designs 1 to 3, 400 samples are generated, and 1000 for design 4. On each of them, a tree
classifier is selected using the growing/pruning strategy, where parameter
$\alpha_n$ is selected by 10-fold cross-validation. Different values of
parameters $n$ ($n=50,100,200$) and noise
($q=0.1, 0.2, 0.3$ in design 1, $\sigma^2=0.5, 1, 2$ in designs 2 and 3, and $\sigma^2=0.2$
in design 4) are used. The number of variables considered to build the
classifiers grows from $p=30$ to
$p=10^3$.\\

\begin{figure}[hbt!]
\begin{center}
\begin{tabular}{cc}
\includegraphics[scale=0.4]{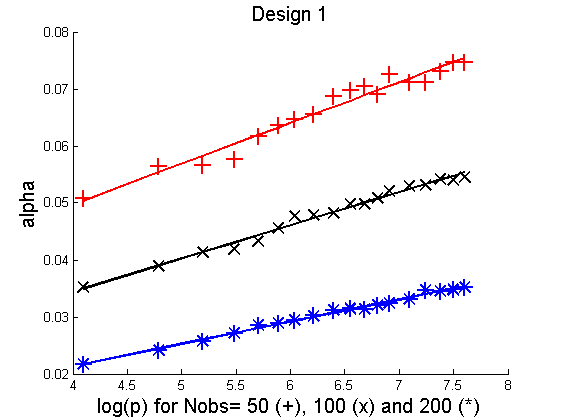}&
\includegraphics[scale=0.4]{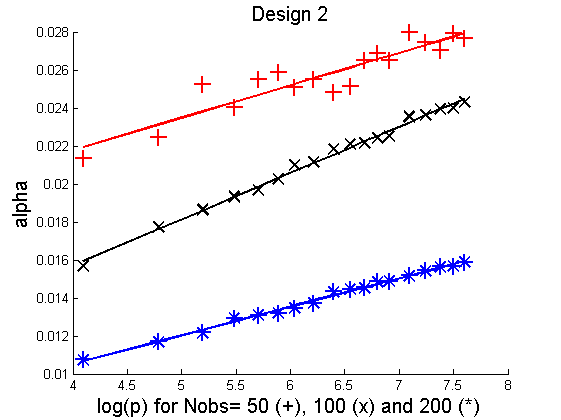}\\
\includegraphics[scale=0.4]{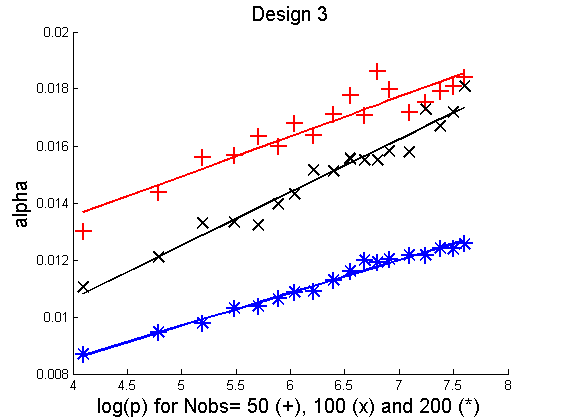}&
\includegraphics[scale=0.4]{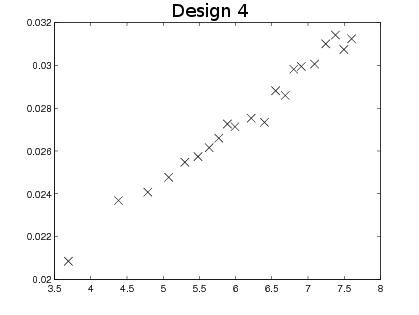}\\
\end{tabular}
\end{center}
\caption{Average value of $\alpha_n$ with respect to  $\log p$, for  $n=50$ (+), $n=100$ (x) and $n=200$ (*). Data are
  simulated from design 1 with $q=0.3$ (Top Left), design 2 (Top Right),
  design 3 (Bottom Left) with $\sigma^2=2$. For design 4 (Bottom Right) the average $\alpha_n$ is obtained over 1000 samples, for $n=100$.}  \label{GraphePenaCART}
\end{figure}

Figure \ref{GraphePenaCART} displays the average value (on 400 simulations) of $\alpha_n$ versus the log-number of
variables for the different designs. Parameter $\alpha_n$
decreases with respect to $n$, and the relationship between the selected
$\alpha_n$ and $\log p$ is linear. These behaviors are observed
whatever the level of noise (not shown) and whatever the design. This confirms
that variable selection is taken into account by the pruning procedure of the
CART algorithm through the choice of $\alpha_n$. This also suggests that the
penalty function proposed in \eqref{Eq : penalite} is relevant regarding
its dependency on $\log p$.

\section{Discussion} \label{Sec : discussion}

As stated in the Introduction, most previous results are related to the pruning step of the CART algorithm rather than considering the general optimization problem \eqref{Equ:PenalizedCriterionMinimization}. For instance, in \cite{Gey12} and \cite{Nob02}, risk bounds are obtained for the collection of CART pruned subtrees, which itself depends on the data at hand: the collection of models includes classes
$\mathcal{C}_{0},...,\mathcal{C}_{{K-1}},\mathcal{C}_{{K}}$ of tree
classifiers built on the maximal tree $T_{K}$, obtained from the training set, and its subtrees
$T_0\preccurlyeq...\preccurlyeq T_{K-1}$. Thus the conditional risk bounds provided in
previous articles only guarantee that the risk of the candidate is
at most of the order of the risk of class $\mathcal{C}_{k^*}$
corresponding to the best subtree $T_{k^*}$. While this exactly
describes the process of the CART algorithm, the guarantee may be
poor if the best subtree of the collection is far from the best tree
among all possible trees. Conversely, the approach presented here
guarantees that the risk bound for the selected tree classifier is
comparable to the risk of the class corresponding to the optimal
tree (among all possible trees). \\

Proposition \ref{Proposition : ModelSelection} generalizes the results obtained in \cite{ScottNowak06} in two ways. First, Scott and Nowak considered the particular case where the tree classifiers are constructed on a fixed dyadic grid. In dyadic trees, the choice of the threshold at each internal node is deterministic, instead of being optimally tuned on the training set. This optimization is taken into account in the results presented here. Second, as recalled in Section \ref{Sec
: Context}, without any margin assumption, the penalty functions obtained in \cite{ScottNowak06} are naturally proportional to
the square root of the tree size over $n$. A $\sqrt{\log{p}}$ factor also appears in the resulting penalties. In comparison, the results presented here exhibit a range of penalty function from square root to linear depending on the margin assumption. If {\bf MA(2)} is satisfied, this validates the form of the penalty implemented in the CART algorithm. If {\bf MA(1)} is satisfied, it leads to better convergence rates for the risk bound. \\

Whenever margin assumption {\bf MA(1)} is satisfied, the penalty suggested in Proposition \ref{Proposition : ModelSelection} is sublinear. In this case the heuristic approach of the CART algorithm can still be employed to obtain an approximate version of $\hat{f}$. Indeed, as proved in \cite{Scott05}, pruning with subadditive penalties
produces sequences of pruned subtrees included in the sequence obtained through pruning with a linear penalty. This means that one can obtain an approximate optimizer of criterion \eqref{Eq : critere}, to the condition that the margin parameter is known. \\

The theoretical form of the penalty term \eqref{Eq : penalite}
derived in Proposition \ref{Proposition : ModelSelection} is of
practical interest. First, it shows that sequential selection
algorithms, such as stepwise or backward variable selection methods,
can be easily studied in the model selection framework where the
selection is supposed to be exhaustive. In the particular case of
tree classification, the simulation study confirms that the penalty
derived under the hypothesis of exhaustive variable selection is the
one that is used in practice by the CART algorithm, that proceeds as
a forward variable selection process. Second, it provides an interesting insight into the CART variable selection process. Indeed, the definition of the classes comes from the fact that a single variable may appear at different nodes, a specificity that changes the classical way of taking into account variable selection in the penalty term: in trees the variable list is ordered (the first variable of the list is associated with the first node) and a
variable may be associated with several nodes. Therefore the classical $\binom{p}{k-1}$ term that appears in penalties in \cite{BirMas01} or \cite{MarRobDau07} (i.e. the number of samplings without replacements and unordered sample) is replaced with
$p^{k-1}$ (i.e. the number of sampling with replacements and ordered sample).\\ 

In \cite{Kol06}, Koltchinskii provides a synthesis of oracle
inequalities in classification. In particular, the author considers
margin assumptions more general than the margin assumption {\bf
MA(1)} given in \cite{MamTsy99}. The in-probability upper bounds for
the loss $l(f^*,\ef)$ given in Propositions \ref{Proposition :
  Control1class} and \ref{Proposition : ModelSelection} can be
straightforwardly generalized using Koltchinskii's margin
definition. This would lead to improved in-probability upper bounds
for the loss $l(f^*,\ef)$, similar to the one given in Theorem 6 of
\cite{Kol06}. However, unlike hypothesis {\bf
  MA(1)} considered here, it would not permit one to obtain explicit rates of
convergence for the risk. Importantly, using a more general margin
assumption would provide no improvement concerning the embedded
selection aspect that we investigated here. From this aspect the
results obtained are tight, as illustrated by the simulation study.

\section*{Acknowledgements}

We would like to thank Sylvain Arlot for helpful discussions and useful
advice.

\section{Proofs} \label{Sec : proofs}

\subsection{Preliminary results} \label{subsec : prelimresults}
We provide two lemmas regarding the Vapnik entropy and the cardinality of tree class collections. \\
Note $H_{c\ell}$ the Vapnik-Chervonenkis log-entropy of class $\mathcal{C}_{c\ell}$:
\begin{eqnarray*}
H_{c\ell}=\log|\{A(f)\cap\{X_1,\ldots,X_n\}, \
f\in \mathcal{C}_{c\ell}\}| ,
\end{eqnarray*}
where $A(f)=\{x \in \mathcal{X}: f(x)=1\}$.
\begin{lemma} \label{Lemma : Entropie}
For a tree class $\mathcal{C}_{c\ell}$, one has
\begin{eqnarray*}
E(H_{c\ell}) \leq |T_{c\ell}|\log(2n)
\end{eqnarray*}
\end{lemma}
\ni This is obtained from lemma (2) in \cite{GeyNed05}. For a tree with $|T_{c\ell}|$ leaves, there are $|T_{c\ell}|-1$ nodes for which the thresholds have to be estimated, leading to at most $n$ ways to split the training sample. The possible number of splittings is bounded by $n^{|T_{c\ell}|-1}$. A given splitting shatters the sample into $|T_{c\ell}|$ subsamples, and each of these subsamples receive label 0 or 1. There are $2^{|T_{c\ell}|}$ ways to label the subsamples, hence
\begin{eqnarray*}
H_{c\ell} &<& \log\left(n^{|T_{c\ell}|-1}\times 2^{|T_{c\ell}|}\right)\\
 &<& |T_{c\ell}|\log(2n) \ \ .
\end{eqnarray*}
Taking the expectation leads to the result.

\begin{lemma}\label{Lemma : NbreModeleParClasse}
The number of classes of trees of size $k$ is $$p^{k-1}N_{(k)} , \ \text{ with } \ \ N_{(k)}=\frac{1}{k}\binom{2k-2}{k-1} \ \ .$$
\end{lemma}
\ni First note that counting the number of classes amounts to counting the number of trees. A tree $T_{c\ell}$ is defined by a configuration $c$ combined with a variable list $\ell$.
The total number of tree configurations of size $k$ is given by the Catalan number $N_{(k)}$.
The total number of lists of $k-1$ variables is $p^{k-1}$, because at each node we have to choose between the $p$ available variables. Combined with the total number of tree configurations, this leads to the proposed lemma.

\paragraph{Remark} In contrast with the classical variable selection framework, in trees the variable list is ordered (the first variable of the
list is associated with the first node) and a
variable may be associated with several nodes. Therefore the classical $\binom{p}{k-1}$ term that appears in penalties in \cite{BirMas01} or \cite{MarRobDau07} (i.e. the number of samplings without replacements and unordered sample) is replaced with
$p^{k-1}$ (i.e. the number of sampling with replacements and ordered sample).

\subsection{Proof of Proposition \ref{Proposition :
    Control1class}} \label{subsec : proofcontrol1class}

A classical way to bound $l(f^*,\widehat{f}_{c\ell})$ is to use the following
decomposition: $$l(f^*,\widehat{f}_{c\ell}) =
l(f^*,\overline{f}_{c\ell})+\P\widehat{f}_{c\ell}-\P\overline{f}_{c\ell},$$
and then to upper bound the variance term
$\P\widehat{f}_{c\ell}-\P\overline{f}_{c\ell}$. In the case where class $\mathcal{C}_{c\ell}$ is finite, an upper bound can be obtained by using Bernstein inequality, as developped in \cite{Lecue07} for instance. In our setting, because there may be (at least) one continuous coordinate (i.e. one continuous variable), classes $\mathcal{C}_{c\ell}$ are not finite. In this case, the upper bounding can be done using Theorem 2 from \cite{Kol06}, which can be restated for our purpose as follows:
\begin{th-a}[Koltchinskii, 2006] \label{Theorem : Kolt}
If there exists a nondecreasing strictly concave function $\psi_{c\ell} : \mathbb{R}_+\rightarrow\mathbb{R}_+$ such that
with probability at least $1-e^{-t_{c\ell}}$
$$ \sup_{f,g \in
\mathcal{C}_{c\ell}(\d)}|(\P_n-\P)(f-g)| \leq\psi_{c\ell}(\delta) \ \  ,$$
and if $\psi^\sharp_{c\ell}$ is defined as
$$\psi^\sharp_{c\ell}(\varepsilon) = \inf\{\delta>0 \text{ s.t. }
\underset{\sigma\geq\delta}{\sup}\frac{\psi(\sigma)}{\sigma}\leq\varepsilon\}
\ \ ,$$ then for all $\d \geq \psi^\sharp_{c\ell}(1/q)$
$$\P\left[\P\widehat{f}_{c\ell}-\P\overline{f}_{c\ell}>\delta\right]\leq
e^{-t_{c\ell}} \ \ .$$\\
\end{th-a}

In order to use Theorem \ref{Theorem : Kolt}, we need to provide an explicit
expression for $\psi_{c\ell}$. To proceed, we start from the following probabilistic upper bound given in
\cite{Kol06} and derived from Talagrand's inequality for bounded processes (see
\cite{Bou02} for more details):
\begin{eqnarray}
\sup_{f,g \in
\mathcal{C}_{c\ell}(\d)}|(\P_n-\P)(f-g)| \leq
2\left(\E\left[\sup_{f,g\in
\mathcal{C}_{c\ell}(\d)}|(\P_n-\P)(f-g)|\right]+D(\mathcal{C}_{c\ell}(\d))\sqrt{\frac{t_{c\ell}}{n}}+\frac{t_{c\ell}}{n}\right)\label{Equ:
BorneSupKolch}
\end{eqnarray}
with probability larger than $1- e^{-t_{c\ell}}$, where $$ \Ccld = \{f\in
\Ccl \text{ s.t. } Pf -  P\overline{f}_{c\ell} \leq\delta\}$$ and
\begin{eqnarray*}
D(\Ccld)&=&\sup_{f,g\in \Ccld}\sqrt{\E((f-g)^2)}\\
&=&\sup_{f,g\in \Ccld}d(f,g)
\end{eqnarray*}
This last term can be upper-bounded in expression \eqref{Equ:
BorneSupKolch} using the margin assumption {\bf MA(1)} described by \eqref{marginTsy}:
\begin{eqnarray*}
d^2(f,f^*) & \leqslant & C_\kappa l(f^*,f)^{\frac{1}{\kappa}}
\end{eqnarray*}
where $C_\kappa = (\kappa-1)^{\frac{1}{\kappa}} C^{\frac{\kappa-1}{\kappa}}\displaystyle{\frac{\kappa}{\kappa-1}}$. Hence
\begin{eqnarray}
d(f,g) & \leqslant &
2\sqrt{C_\kappa}\left(l(f^*,\overline{f}_{c\ell})^{\frac{1}{2\kappa}} +
\delta^{\frac{1}{2\kappa}}\right)\nonumber\\
\Rightarrow D(\Ccld)& \leqslant &
2\sqrt{C_\kappa}\left(l(f^*,\overline{f}_{c\ell})^{\frac{1}{2\kappa}} +
\delta^{\frac{1}{2\kappa}}\right) = D \ \ . \label{Equ :
MajorationLabel}
\end{eqnarray}
Now because
\begin{eqnarray*}
\E\left[\sup_{f,g\in
\mathcal{C}_{c\ell}(\d)}|(\P_n-\P)(f-g)|\right] \leq
\E\left[\sup_{d(f,g)\leq D}|(\P_n-\P)(f-g)| \right]
\end{eqnarray*}
we can use the result of \cite{Mas00} (p295) to obtain
\begin{eqnarray}
\E\left[\sup_{f,g\in
\mathcal{C}_{c\ell}(\d)}|(\P_n-\P)(f-g)|\right] \leq
24D\sqrt{\frac{E[H_{c\ell}]}{n}} \ \ , \label{Equ : MajorationBiais1}
\end{eqnarray}
where $H_{c\ell}$ is the Vapnik-Chervonenkis log-entropy of
$\mathcal{C}_{c\ell}$. Combining (\ref{Equ : MajorationLabel}) and (\ref{Equ :
MajorationBiais1}), then using lemma 5 of \cite{TsyGee05}, we
obtain for all $\alpha\in]0,1[$
\begin{eqnarray}
\sup_{f,g \in
\mathcal{C}_{c\ell}(\d)}|(\P_n-\P)(f-g)| &\leq&
2\left[2\sqrt{C_\kappa}\left(24\sqrt{\frac{E[H_{c\ell}]}{n}} +
\sqrt{\frac{t_{c\ell}}{n}}\right)\left(l(f^*,\overline{f}_{c\ell})^{\frac{1}{2\kappa}}
+ \delta^{\frac{1}{2\kappa}}\right)+\frac{t_{c\ell}}{n}\right]\nonumber\\
&\leq& 4\sqrt{C_\kappa}\left(24\sqrt{\frac{E[H_{c\ell}]}{n}} +
\sqrt{\frac{t_{c\ell}}{n}}\right)\delta^{\frac{1}{2\kappa}}\nonumber\\
& &+ 2\frac{t_{c\ell}}{n} + \alpha l(f^*,\overline{f}_{c\ell}) +
\beta_{\kappa,\alpha}\left(\frac{E[H_{c\ell}]}{n}\right)^{\frac{\kappa}{2\kappa-1}}+
\frac{\beta_{\kappa,\alpha}}{24}\left(\frac{t_{c\ell}}{n}\right)^{\frac{\kappa}{2\kappa-1}}
\ \ .
\nonumber
\end{eqnarray}

\noindent In the present framework, we then have
\begin{eqnarray*}
\psi_{c\ell}(\delta) &=& 4\sqrt{C_\kappa}\left(24\sqrt{\frac{E[H_{c\ell}]}{n}}
+ \sqrt{\frac{t_{c\ell}}{n}}\right)\delta^{\frac{1}{2\kappa}} +
2\frac{t_{c\ell}}{n} + \alpha l(f^*,\overline{f}_{c\ell}) +
\beta_{\kappa,\alpha}
\left(\frac{E[H_{c\ell}]}{n}\right)^{\frac{\kappa}{2\kappa-1}}+
\frac{\beta_{\kappa,\alpha}}{24}\left(\frac{t_{c\ell}}{n}\right)^{\frac{\kappa}{2\kappa-1}}\\
&=& \psi_1(\delta) + \psi_2(\delta) + K
\end{eqnarray*}
where
\begin{eqnarray*}
\psi_1(\delta) &=&
96\sqrt{C_\kappa}\sqrt{\frac{E[H_{c\ell}]}{n}}
\delta^{\frac{1}{2\kappa}} \\
\psi_2(\delta) &=&
4\sqrt{C_\kappa}\sqrt{\frac{t_{c\ell}}{n}}\delta^{\frac{1}{2\kappa}} \\
\text{ and } K &=& 2\frac{t_{c\ell}}{n} + \alpha
l(f^*,\overline{f}_{c\ell}) + \beta_{\kappa,\alpha}
\left(\frac{E[H_{c\ell}]}{n}\right)^{\frac{\kappa}{2\kappa-1}}+
\frac{\beta_{\kappa,\alpha}}{24}\left(\frac{t_{c\ell}}{n}\right)^{\frac{\kappa}{2\kappa-1}}
\end{eqnarray*}
Moreover, $\psi_{c\ell}^\sharp(\varepsilon) \leq \psi^\sharp_1(\varepsilon/3) + \psi^\sharp_1(\varepsilon/3) + \frac{3K}{\varepsilon}$, and $\psi^\sharp_1$ and $\psi^\sharp_2$ can be determined using the following characterization (available for all strictly concave functions $\psi$):
\begin{eqnarray*}
\forall \varepsilon> 0, \ \psi\left(\psi^\sharp(\varepsilon)\right) = \psi^\sharp(\varepsilon)\varepsilon \ \ .
\end{eqnarray*}
Solving this last equation for the particular form of functions $\psi_1$ and $\psi_2$, we obtain
\begin{eqnarray*}
\psi^\sharp_{c\ell}(\varepsilon) &\leq& \left(\frac{288\sqrt{C_\kappa}\sqrt{E[H_{c\ell}]}}{\varepsilon\sqrt{n}}\right)^\frac{2\kappa}{2\kappa-1}
+ \left(\frac{12\sqrt{C_\kappa}\sqrt{t_{c\ell}}}{\varepsilon\sqrt{n}}\right)^\frac{2\kappa}{2\kappa-1}\\
&&+ \left(2\frac{t_{c\ell}}{n} + \alpha l(f^*,\overline{f}_{c\ell}) + \beta_{\kappa,\alpha}
\left(\frac{E[H_{c\ell}]}{n}\right)^{\frac{\kappa}{2\kappa-1}}+
\frac{\beta_{\kappa,\alpha}}{24}\left(\frac{t_{c\ell}}{n}\right)^{\frac{\kappa}{2\kappa-1}}\right)\frac{3}{\varepsilon}
\end{eqnarray*}
Taking $\varepsilon = 1/q$ one has with probability larger than $1-e^{-t_{c\ell}}$
\begin{eqnarray*}
\P\widehat{f}_{c\ell}-\P\overline{f}_{c\ell} &\leq&
\left(\left(q288\sqrt{C_\kappa}\right)^\frac{2\kappa}{2\kappa-1} + 3q \beta_{\kappa,\alpha}
\right) \left(\frac{E[H_{c\ell}]}{n}\right)^\frac{\kappa}{2\kappa-1}
+\left(\left(q12\sqrt{C_\kappa}\right)^\frac{2\kappa}{2\kappa-1}
+ \frac{3q\beta_{\kappa,\alpha}}{24}\right)
\left(\frac{t_{c\ell}}{n}\right)^\frac{\kappa}{2\kappa-1}\\
&&+6q\frac{t_{c\ell}}{n} + 3q\alpha l(f^*,\overline{f}_{c\ell}) \ \ .
\end{eqnarray*}
Using Lemma \ref{Lemma : Entropie} and rescaling $\alpha$ properly, this leads to
\begin{eqnarray}
l(f^*,\widehat{f}_{c\ell}) \leq (1+\alpha)l(f^*,\overline{f}_{c\ell}) + K_{\alpha,\kappa,q}^1\left(\frac{|T_{c\ell}|\log(2n)}{n}\right)^\frac{\kappa}{2\kappa-1} + K_{\alpha,\kappa,q}^2\left(\frac{t_{c\ell}}{n}\right)^\frac{\kappa}{2\kappa-1} + K_q\frac{t_{c\ell}}{n} \ \ .
\label{Equ : LesConstantes}
\end{eqnarray}
Renaming $K_{\alpha,\kappa,q}^1=K_1$, $K_{\alpha,\kappa,q}^2=K_2$ and $K_q=K$ leads to the first expression in Proposition \ref{Proposition : Control1class}. The risk bound follows by integration.

\subsection{Proof of Proposition \ref{Proposition :
    ModelSelection}} \label{subsec : proofmodelselection}
We first choose the weights $t_{c\ell}=x_{c\ell}+x$ associated with classes $\mathcal{C}_{c\ell}$ such that $x_{c\ell}$ and $x$ are positive and
\begin{eqnarray*}
\sum_{c,\ell}e^{-x_{c\ell}} = \Sigma < +\infty \ \ .
\end{eqnarray*}
The exact form of the weights will be chosen later. Furthermore, we will use lemma 4 of \cite{Kol06}, reformulated here for our purpose:
\begin{lemma-a}[Koltchinskii, 2006] \label{Lemma : MajorationBiais}
Consider a class $\mathcal{C}_{c\ell}$ and assume that {\bf MA(1)} is satisfied. For all $t_{c\ell}>0$ and $\alpha \in ]0,2/5[$, with probability at least $1-2e^{-t_{c\ell}}$, one has
\begin{eqnarray}
\P_n\overline{f}_{c\ell}-\P_nf^* \leq (1+\alpha)(\P \overline{f}_{c\ell} - \P f^*) + K_\alpha\left(\frac{t_{c\ell}}{n}\right)^\frac{\kappa}{2\kappa-1} + \frac{t_{c\ell}}{n} \label{Equ : MajorationBiaisEmpirique}
\end{eqnarray}
and
\begin{eqnarray}
\P \overline{f}_{c\ell}-\P f^* \leq \left(1-\frac{5}{2}\alpha\right)^{-1}\left( \P_n\widehat{f}_{c\ell}-\P_nf^*
+ \frac{3}{2}K_1\left(\frac{|T_{c\ell}|\log(2n)}{n}\right)^\frac{\kappa}{2\kappa-1} + 3K_2\left(\frac{t_{c\ell}}{n}\right)^\frac{\kappa}{2\kappa-1} + 3K\frac{t_{c\ell}}{n}
\right) \label{Equ : MajorationBiais}
\end{eqnarray}
with the same notations as above.
\end{lemma-a}
We start the proof from the result obtained in Proposition \ref{Proposition : Control1class}. Combining equation \eqref{PropControl1Class : Traject} of Proposition \ref{Proposition : Control1class} and a classical union bound argument, one has with probability larger than $1-\Sigma e^{-x}$
\begin{eqnarray*}
l(f^*,\tilde{f}) \leq (1+\alpha)l(f^*,\overline{f}_{\widehat{c\ell}}) + K_1\left(\frac{|T_{\widehat{c\ell}}|\log(2n)}{n}\right)^\frac{\kappa}{2\kappa-1} + K_2\left(\frac{x_{\widehat{c\ell}}+x}{n}\right)^\frac{\kappa}{2\kappa-1} + K\frac{x_{\widehat{c\ell}}+x}{n} \ \ ,
\end{eqnarray*}
where $\alpha \in ]0,2/5[$. We now use equation \eqref{Equ : MajorationBiais} from Lemma \ref{Lemma : MajorationBiais} to obtain with probability larger than $1-3\Sigma e^{-x}$
\begin{eqnarray*}
l(f^*,\tilde{f}) &\leq& \frac{(1+\alpha)}{1-\frac{5\alpha}{2}}\left( \P_n\widehat{f}_{\widehat{c\ell}}-\P_nf^* + \frac{5K_1}{2}\left(\frac{|T_{\widehat{c\ell}}|\log(2n)}{n}\right)^\frac{\kappa}{2\kappa-1} + 4K_2\left(\frac{x_{\widehat{c\ell}}+x}{n}\right)^\frac{\kappa}{2\kappa-1} +
4K\frac{x_{\widehat{c\ell}}+x}{n}\right)\\
&\leq& \frac{(1+\alpha)}{1-\frac{5\alpha}{2}}\left( \P_n\widehat{f}_{\widehat{c\ell}}-\P_nf^* + \frac{5K_1}{2}\left(\frac{|T_{\widehat{c\ell}}|\log(2n)}{n}\right)^\frac{\kappa}{2\kappa-1} + 4K_2\left(\frac{x_{\widehat{c\ell}}}{n}\right)^\frac{\kappa}{2\kappa-1} +
4K\frac{x_{\widehat{c\ell}}}{n}\right) \\
&&+ \frac{(1+\alpha)}{1-\frac{5\alpha}{2}}\left( 4K_2\left(\frac{x}{n}\right)^\frac{\kappa}{2\kappa-1} +
4K\frac{x}{n}\right)
\end{eqnarray*}
In the context of variable selection, one has to choose the weights such that
\begin{eqnarray*}
\sum_{c,\ell}e^{-x_{c\ell}} < +\infty \Rightarrow \sum_{k}\sum_{\mathcal{C}_{c\ell} \text{ s.t. } |T_{c\ell}|=k}e^{-x_{c\ell}}< +\infty \ \ .
\end{eqnarray*}
Giving equal weights $x_k$ to classes of same complexity $k$ (i.e. classes $\mathcal{C}_{c\ell}$ and $\mathcal{C}_{c'\ell'}$ such that $|T_{c\ell}|=|T_{c'\ell'}|=k$), one obtains from Lemma \ref{Lemma : NbreModeleParClasse}:
\begin{eqnarray*}
\sum_{k}\sum_{\mathcal{C}_{c\ell} \text{ s.t. } |T_{c\ell}|=k}e^{-x_{c\ell}} &=& \sum_{k} p^{k-1}\frac{1}{k}\binom{2k-2}{k-1} e^{-x_k} \\
&\leq& \sum_{k} \frac{(4p)^k}{k} e^{-x_k} \ \ .
\end{eqnarray*}
The choice $x_{c\ell} = x_{|T_{c\ell}|} = \lambda|T_{c\ell}|\log(p)$ with
$\lambda>3$ ensures that the sum is finite. Hence,
\begin{eqnarray*}
l(f^*,\tilde{f}) &\leq& \frac{(1+\alpha)}{1-\frac{5\alpha}{2}}\left( \P_n\widehat{f}_{\widehat{c\ell}}-\P_nf^* + \frac{5K_1}{2}\left(\frac{|T_{\widehat{c\ell}}|\log(2n)}{n}\right)^\frac{\kappa}{2\kappa-1} +
4K_2\left(\frac{\lambda|T_{\widehat{c\ell}}|\log(p)}{n}\right)^\frac{\kappa}{2\kappa-1} \right. \\
&&\left.+4K\frac{\lambda|T_{\widehat{c\ell}}|\log(p)}{n}\right)+ \frac{(1+\alpha)}{1-\frac{5\alpha}{2}}\left( 4K_2\left(\frac{x}{n}\right)^\frac{\kappa}{2\kappa-1} +
4K\frac{x}{n}\right) \\
&\leq&  \frac{(1+\alpha)}{1-\frac{5\alpha}{2}}\left( \P_n\widehat{f}_{\widehat{c\ell}}-\P_nf^* + C'_\kappa\left(\frac{|T_{\widehat{c\ell}}|\log(2n)}{n}\right)^\frac{\kappa}{2\kappa-1} + C''_\kappa\left(\frac{|T_{\widehat{c\ell}}|\log(p)}{n}\right)^\frac{\kappa}{2\kappa-1} +
C'''_\kappa\left(\frac{|T_{\widehat{c\ell}}|\log(p)}{n}\right) \right)\\
&&+
\frac{(1+\alpha)}{1-\frac{5\alpha}{2}}\left( 4K_2\left(\frac{x}{n}\right)^\frac{\kappa}{2\kappa-1} +
4K\frac{x}{n}\right) \ \ ,
\end{eqnarray*}
for a proper choice of constants $C'_\kappa$, $C''_\kappa$, and $C'''_\kappa$. This leads to
\begin{eqnarray*}
l(f^*,\tilde{f}) \leq \frac{(1+\alpha)}{1-\frac{5\alpha}{2}}\underset{c,\ell}{\inf}\left(\P_n\widehat{f}_{c\ell}-\P_nf^*+pen(c,\ell) \right) + \frac{(1+\alpha)}{1-\frac{5\alpha}{2}}\left( 4K_2\left(\frac{x}{n}\right)^\frac{\kappa}{2\kappa-1} +
4K\frac{x}{n}\right) \ \ .
\end{eqnarray*}
Since $\P_n\widehat{f}_{c\ell}-\P_nf^* \leq \P_n\overline{f}_{c\ell}-\P_nf^*$ ( by definition of $\widehat{f}_{c\ell}$), this last expression can be upper bounded (with probability larger than $1-3\Sigma e^{-x}$) thanks to equation \eqref{Equ : MajorationBiaisEmpirique} of Lemma \ref{Lemma : MajorationBiais}:
\begin{eqnarray*}
l(f^*,\tilde{f}) &\leq& \frac{(1+\alpha)^2}{1-\frac{5\alpha}{2}}\underset{c,\ell}{\inf}\left(\P \overline{f}_{c\ell} - \P f^* + K_\alpha\left(\frac{x_{c\ell}}{n}\right)^\frac{\kappa}{2\kappa-1}+\left(\frac{x_{c\ell}}{n}\right)
+K_\alpha\left(\frac{x}{n}\right)^\frac{\kappa}{2\kappa-1}+\left(\frac{x}{n}\right) +pen(c,\ell) \right)\\
&&+ \frac{(1+\alpha)}{1-\frac{5\alpha}{2}}\left( 4K_2\left(\frac{x}{n}\right)^\frac{\kappa}{2\kappa-1} +
4K\frac{x}{n}\right) \\
&\leq& \frac{2(1+\alpha)^2}{1-\frac{5\alpha}{2}}\underset{c,\ell}{\inf}\left(\P \overline{f}_{c\ell} - \P f^* +
pen(c,\ell) \right)+ \frac{2(1+\alpha)^2}{1-\frac{5\alpha}{2}}\left( 4K_2\left(\frac{x}{n}\right)^\frac{\kappa}{2\kappa-1} +
4K\frac{x}{n}\right) \\
&\leq& C_1'\underset{c,\ell}{\inf}\left(\P \overline{f}_{c\ell} - \P f^* +
pen(c,\ell) \right)+ C'_2\left( \left(\frac{x}{n}\right)^\frac{\kappa}{2\kappa-1} +
\frac{x}{n}\right) \ \ .
\end{eqnarray*}
The last inequality corresponds to the first equation of Proposition \ref{Proposition : ModelSelection}. The risk bound follows by integration.

\bibliographystyle{plain}
\bibliography{pruning}

\end{document}